\documentclass{amsart}
\usepackage{amsmath,amssymb,hyperref}
\usepackage[margin=1.5in]{geometry}
\usepackage[vcentermath]{youngtab}
\def\CC{\mathbb{C}} 
\def\S{\mathfrak{S}}
 
\def\<{\langle}
\def\>{\rangle}

\DeclareMathOperator{\End}{End}

\DeclareMathOperator{\GL}{GL}

\theoremstyle{definition}
\newtheorem{thm}{Theorem}
\newtheorem{cor}[thm]{Corollary}
\newtheorem{prop}[thm]{Proposition}

\begin{document}
\title{A Short Proof of Gamas's Theorem}
\author{Andrew Berget}
\address{School of Mathematics,
  University of Minnesota,
  Minneapolis, MN 55455}
\email{berget@math.umn.edu}

\thanks{The author is grateful to Victor Reiner for listening to the
  many incarnations of this proof. His support from NSF grant DMS
  0601010 was also greatly appreciated.\\
  {\url{http://dx.doi.org/10.1016/j.laa.2008.09.027}}}

\begin{abstract}
  If $\chi^\lambda$ is the irreducible character of $\mathfrak{S}_n$
  corresponding to the partition $\lambda$ of $n$ then we may
  symmetrize a tensor $v_1 \otimes \cdots \otimes v_n$ by
  $\chi^\lambda$.  Gamas's theorem states that the result is not zero
  if and only if we can partition the set $\{v_i\}$ into linearly
  independent sets whose sizes are the parts of the transpose of
  $\lambda$. We give a short and self-contained proof of this fact.
\end{abstract}
\maketitle
\section{Introduction}
Let $\lambda$ be a partition of a positive integer $n$ and let
$\chi^\lambda$ be the irreducible character of the symmetric group
$\S_n$ corresponding to $\lambda$. There is a right action of $\S_n$
on $V^{\otimes n}$, where $V$ is a finite-dimensional complex vector
space, by permuting tensor positions. Let $T_\lambda$ be the
endomorphism of $V^{\otimes n}$ given by
\[
(v_1 \otimes \dots \otimes v_n)
T_\lambda=\frac{\chi^\lambda(1)}{n!}\sum_{\sigma \in \S_n}
\chi_\lambda(\sigma) v_{\sigma(1)} \otimes \dots \otimes
v_{\sigma(n)}.
\]
Our goal is to prove the following result of Carlos Gamas \cite{gamas}.
\begin{thm}[Gamas's Theorem]
  Let $v_1, \dots, v_n$ be vectors in $V$. Then
  \[
  (v_1 \otimes \cdots \otimes v_n)T_\lambda \neq 0
  \]
  if and only if it is possible to partition the set $\{v_i\}$ into
  linearly independent sets whose sizes are the parts of the transpose
  of $\lambda$.
\end{thm}
If $\{v_1,\dots,v_n\}$ is a collection of vectors satisfying the
condition of the theorem we will say that it satisfies ``Gamas's
Condition for $\lambda$''. The theorem is a generalization of the well
known fact that the exterior product of a set of vectors is nonzero if
and only the set of vectors is linearly independent.

In addition to Gamas's proof of this result there was a second one
given by Pate in 1990 \cite{pate} using results he obtained in
\cite{pate2}. The benefit of our proof is that it is self-contained
and short. It relies on standard facts from the representation theory
of $\GL(V)$, namely, Schur-Weyl duality and the Pieri Rule. We refer
to Fulton and Harris's book \cite{fulton-harris} for the needed
background and notation.

\section{Preliminaries and Proof}
Let $V$ be a finite dimensional complex vector space. The general
linear group $\GL(V)$ acts diagonally on $V^{\otimes n}$. Let $w \in
V^{\otimes n}$ be any tensor. Define $G(w)$ to be the $\GL(V)$-module
spanned by
\[
\GL(V)w = \{ g\cdot w : g \in \GL(V)\}.
\]
We are interested in which irreducible $\GL(V)$-modules appear in
$G(w)$. Since $G(w) \subset V^{\otimes n}$ is a polynomial
representation, the isomorphism type of the irreducible
$\GL(V)$-modules which can appear in $G(w)$ are indexed by partitions
$\lambda$ of $n$ with at most $\dim V$ parts.  If $\lambda$ is such a
partition, we will say that \textit{$\lambda$ appears in $G(w)$} if
this module has a highest weight vector of weight $\lambda$ (see
\cite[Chapter~15]{fulton-harris}). We will write $\ell(\lambda)$ for
the number of parts of $\lambda$.
\begin{prop}\label{prop}
  If $\lambda$ is a partition of $n$, then $\lambda$ appears in $G(w)$
  if and only if $w T_\lambda \neq 0$.
\end{prop}
\textit{Proof.} Note that $\lambda$ appears in $G(w)$ if and only if
the projection of $G(w)$ onto its $\lambda$-th isotypic component is
not zero. By Schur-Weyl duality (see \cite[Lemma~6.22]{fulton-harris})
$T_\lambda$ is this projector, since it is the projector of $\CC\S_n$
onto its $\lambda$-th isotypic component. Since $T_\lambda$ commutes
with the $\GL(V)$ action, this isotypic component is zero if and only
if $G(w T_\lambda) = 0$, which happens if and only if $w T_\lambda
=0$.  \qed

The following corollary is immediate from Proposition~\ref{prop}.
\begin{cor}\label{cor:diminvariance}
  Suppose that $W$ is a subspace of $V$ and $w \in W^{\otimes n}
  \subset V^{\otimes n}$.  The shape $\lambda$ appears in
  $\operatorname{span}\GL(V) w$ if and only if it appears in
  $\operatorname{span}\GL(W) w$.
\end{cor}

\textit{Proof of Gamas's Theorem.} Assume that the vectors
$\{v_1,\dots,v_n\}$ span $V$, as we may by
Corollary~\ref{cor:diminvariance}. Suppose that $\{v_i\}$ satisfy
Gamas's condition for $\lambda$. We prove the result by induction on
$n+\ell(\lambda)$. Write $v^\otimes$ for the tensor $v_1 \otimes \dots
\otimes v_n$.

If $\lambda$ has one part $\chi^\lambda$ is the trivial character and
$ v^\otimes T_\lambda$ is a scalar multiple of the fully symmetrized
tensor $v_1 \cdots v_n$ in $\operatorname{Sym}^n(V)$. This is not zero
since none of the $v_i$ are zero.

If $\ell(\lambda)< \dim V$ let $A \in \End(V)$ be a generic projection
to a subspace $W \subset V$ with dimension equal to the length of
$\lambda$.  Since $A$ is generic, the collection $\{Av_1,\dots,Av_n\}$
still satisfies Gamas's condition for $\lambda$.  It follows by
induction that $\lambda$ appears in the span of $\GL(W)(A\cdot
v^\otimes)$ and hence it also appears in $G(A \cdot v^\otimes)$. Since
$A$ is a limit of elements of $\GL(V)$ we have $G(A\cdot v^\otimes)
\subset G(v^\otimes)$ and hence $\lambda$ appears in $G(v^\otimes)$.

If $\ell(\lambda) =\dim V$ then we consider a Young tableau of shape
$\lambda$ whose columns index independent subsets of the set $v=\{v_1,
\dots ,v_n\}$. Let $B$ be the set of numbers in the first column of
the tableau. The map
\[
b_{B}=\sum_{\sigma \in \S_{B}} (-1)^\sigma \sigma \in
\CC\mathfrak{S}_n = \operatorname{End}_{\GL(V)}(V^{\otimes n})
\]
is a map of $\GL(V)$-modules and hence we have a surjection of
$\GL(V)$-modules $G(v^\otimes) \to G( v^\otimes b_{B})$.  Without loss
of generality, we write $B = \{1,\dots,k\}$ where $k = \dim V$, so
that
\[
G( v^\otimes b_B)= \textstyle\det_V \otimes G(v_{k+1}\otimes \dots
\otimes v_n).
\]
Here $\det_V$ is the one dimensional representation $g \mapsto
\det(g)$ of $\GL(V)$. For example, if $\dim V = 2$ and $B = \{1,2\}$
then
\[
(v_1 \otimes v_2 \otimes v_3 \otimes v_4 \otimes v_5) b_B = (v_1
\otimes v_2 - v_2 \otimes v_1) \otimes v_3 \otimes v_4 \otimes v_5
\]
Since $v_1$ and $v_2$ are a basis of $V$, we see that $g \in \GL(V)$
acts by its determinant on the exterior power $\bigwedge^2 V$ and
hence on $v_1 \wedge v_2 = v_1 \otimes v_2 - v_2 \otimes v_1$.

Denote by $\lambda^-$ the shape obtained from $\lambda$ by removing
the first column. Then $\{v_{k + 1},\dots,v_n\}$ satisfies Gamas's
condition for $\lambda^-$. By induction we know that ${\lambda^-}$
appears in $G( v_{k+1}\otimes \dots \otimes v_n)$. By Pieri's Rule
(see \cite[Equation~6.9]{fulton-harris}) it follows that $\lambda$
appears in
\[
\textstyle \det_V \otimes  G( v_{k+1}\otimes \dots \otimes v_n)
\]
and, hence, in $G(v^\otimes)$ since it appears in its homomorphic
image $G(v^\otimes b_B)$. This completes the more difficult
implication of Gamas's Theorem.

Although our proof of the converse was essentially known to Pate
\cite{pate}, we include it to keep this paper self-contained. We will
need the standard construction of the irreducible $\GL(V)$ and
$\CC\S_n$ modules via Young symmetrizers. To this end let $T$ be a
tableau of shape $\lambda$, $a_T$ its row symmetrizer, and $b_T$ its
column antisymmetrizer. These are given by
\[
\sum_{\sigma \in \operatorname{Row}(T)} \sigma, \quad \sum_{\sigma \in
  \operatorname{Col}(T)} \operatorname{sign}(\sigma)\sigma,
\]
respectively. For example, using cycle notation for permutations in
$\S_n$, if
\[
T = \young(234,15) 
\]
then $b_T = (1-(12))(1-(35))$ while
\[
a_T = \left(1+(23)+(24)+(34)+(234)+(243)\right)(1+(15)).
\]
A product $b_Ta_T$ is called a \textit{Young symmetrizer} and the
right ideal in $\CC\S_n$ generated by a Young symmetrizer is an
irreducible $\CC\S_n$-module with character $\chi^\lambda$ while the
image of $b_T a_T$ on $V^{\otimes n}$ is zero, or irreducible with
highest weight $\lambda$ (see \cite[Chapters 4 and
15]{fulton-harris}). It is clear that $v^\otimes b_T$ is not zero if
and only if the sets of vectors indexed by the columns of the tableau
$T$ are linearly independent.

It follows from Schur-Weyl duality that if $\lambda$ appears in
$G(v^\otimes)$ then there is an element $c \in
\End_{\GL(V)}(V^{\otimes n})= \CC\S_n$ such that $G(v^\otimes c)$
equals the irreducible $\GL(V)$-module $V^{\otimes n}b_T a_T$. It then
follows that the right $\CC\S_n$-module generated by $v^\otimes c$ is
isomorphic to both $c \CC\S_n$ and $b_T a_T \CC\S_n$, in particular
the latter two modules are isomorphic. We conclude that $c$ can be
written as a sum $c = \sum_{\sigma \in \S_n} x_\sigma \sigma b_Ta_T$,
$x_\sigma \in \CC$, and hence one of these terms $x_\sigma \sigma b_T
a_T$ applied to $v^\otimes$ is not zero.  Finally, since $v^\otimes
\sigma b_T$ is not zero Gamas's Condition holds for $\lambda$, the
shape of $T$.\qed

Define a sequence of integers $\rho_i$ by the condition that
\[
\rho_1 + \dots + \rho_k
\]
is the size of the largest union of $k$ linearly independent subsets
of $\{v_i\}$. The sequence $\rho$ is called the \textit{rank
  partition} of $\{v_i\}$ and was introduced by Dias da Silva in
\cite{dds}. In our language, Dias da Silva proved the following
strengthening of Gamas's Theorem.
\begin{thm}[Dias da Silva]
  The partition $\lambda$ appears in $G(v^\otimes)$ if and only if
  $\lambda$ is larger (in dominance order) than the transposed rank
  partition of $\{v_i\}$.
\end{thm}
The extent to which one can further predict the irreducible
$\GL(V)$-decomposition of $G(v^\otimes)$ is the subject of the
author's Ph.D. thesis.

\bibliographystyle{plain}
\bibliography{gamas-thm}

\begin{thebibliography}{1}

\bibitem{dds}
J.A. Dias~da Silva.
\newblock On the {$\mu$}-colorings of a matroid.
\newblock {\em Linear and Multilinear Algebra}, 27(1):25--32, 1990.

\bibitem{fulton-harris}
William Fulton and Joe Harris.
\newblock {\em Representation theory}, volume 129 of {\em Graduate Texts in
  Mathematics}.
\newblock Springer-Verlag, New York, 1991.
\newblock A first course, Readings in Mathematics.

\bibitem{gamas}
Carlos Gamas.
\newblock Conditions for a symmetrized decomposable tensor to be zero.
\newblock {\em Linear Algebra Appl.}, 108:83--119, 1988.

\bibitem{pate}
Thomas~H. Pate.
\newblock Immanants and decomposable tensors that symmetrize to zero.
\newblock {\em Linear and Multilinear Algebra}, 28(3):175--184, 1990.

\bibitem{pate2}
Thomas~H. Pate.
\newblock Partitions, irreducible characters, and inequalities for generalized
  matrix functions.
\newblock {\em Trans. Amer. Math. Soc.}, 325(2):875--894, 1991.

\end{thebibliography}

\end{document}